\documentclass[MSNbibl,number,citesort,seceqn,dvips]{arxbj}
\usepackage{upgreek}
\usepackage{mathrsfs}
\usepackage{mathbh}

\aid{0}
\volume{18}
\issue{4}
\pubyear{2012}
\firstpage{1267}
\lastpage{1283}
\doi{10.3150/11-BEJ371}

\makeatletter
\newcommand{\eqref}[1]{(\ref{#1})}
\newtheorem{prp}{Proposition}
\newremark{remark}{Remark}
\makeatother

\begin{document}
\begin{frontmatter}

\title{Asymptotics for a Bayesian nonparametric estimator of species variety}
\runtitle{Asymptotics for species variety}

\begin{aug}
\author[a,c]{\fnms{Stefano} \snm{Favaro}\thanksref{a,c,e1}\ead[label=e1,mark]{stefano.favaro@unito.it}},
\author[b,c]{\fnms{Antonio} \snm{Lijoi}\corref{}\thanksref{b,c}\ead[label=e2]{lijoi@unipv.it}}
\and
\author[a,c]{\fnms{Igor} \snm{Pr\"unster}\thanksref{a,c,e3}\ead[label=e3,mark]{igor@econ.unito.it}}

\address[a]{Dipartimento di Statistica e Matematica Applicata,
Universit\`a degli Studi di Torino C.so Unione Sovietica 218/bis,
10134 Torino, Italy. \printead{e1,e3}}

\address[b]{Dipartimento di Economia Politica
e Metodi Quantitativi, Universit\`a degli Studi di Pavia, Via San
Felice~5, 27100 Pavia, Italy. \printead{e2}}

\address[c]{Collegio Carlo Alberto, Via Real Collegio 30, 10024
Moncalieri, Italy}

\runauthor{S. Favaro, A. Lijoi and I. Pr\"unster}
\end{aug}

\received{\smonth{8} \syear{2010}}
\revised{\smonth{1} \syear{2011}}

%
\begin{abstract}
In Bayesian nonparametric inference, random discrete probability
measures are
commonly used as priors within hierarchical mixture models for density
estimation
and for inference on the clustering of the data. Recently, it has been
shown that they can also be exploited in species sampling problems:
indeed they are natural tools for modeling the random proportions of
species within a population thus allowing for inference on various
quantities of statistical interest. For applications that involve large
samples, the exact evaluation of the corresponding estimators becomes
impracticable and, therefore, asymptotic approximations are sought. In
the present paper, we study the limiting behaviour of the number of new
species to be observed from further sampling, conditional on observed
data, assuming the observations are exchangeable and directed by a
normalized generalized gamma process prior. Such an asymptotic study
highlights a connection between the normalized generalized gamma
process and the two-parameter Poisson--Dirichlet process that was
previously known only in the unconditional case.
\end{abstract}

%

\end{frontmatter}
%

\section{Introduction}\label{section1}
In species sampling problems, one is interested in the species
composition of a certain population (of plants, animals, genes, etc.)
containing an unknown number of species and only a sample drawn from it
is available. The relevance of such problems in ecology, biology and,
more recently, in genomics and bioinformatics is not surprising. From
an inferential perspective, one is willing to use available data in
order to evaluate some quantities of practical interest. The available
data specifically consist of a so-called \textit{basic sample} of size
$n$, $(X_1, \ldots, X_n)$, which exhibits $K_n \in\{1, \ldots, n\}$
distinct species, $(X_1^*, \ldots, X_{K_n}^*)$, with respective
frequencies $(N_1, \ldots, N_{K_n})$, where clearly $\sum_{i=1}^{K_n}
N_i=n$. Given a basic sample, interest mainly lies in estimating the
number of new species, $K_m^{(n)}:=K_{m+n}-K_n$, to be observed in an
additional sample $(X_{n+1},\ldots,X_{n+m})$ of size~$m$ and not
included among the $X_j^*$'s, $j=1,\ldots,K_n$.

Most of the contributions in the literature that address this issue
rely on a frequentist approach (see \cite{BF93,C05} for reviews)
and only recently an alternative Bayesian nonparametric approach has
been set forth (see, e.g., \cite{Fav09,Lij207,Lij08,gnedin2}).
The latter resorts to a general class of discrete random probability
measures, termed \textit{species sampling models} and introduced by
J.~Pitman in \cite{Pit96a}. Given a nonatomic probability measure
$P_0$ on some complete and separable metric space $\mathbb{X}$,
endowed with
the Borel $\sigma$-field $\mathscr{X}$, a (proper) species sampling
model on
$(\mathbb{X}, \mathscr{X})$ is a random probability measure
\[
\tilde p=\sum_{i\ge1} \tilde p_i \delta_{Y_i},
\]
where $(Y_i)_{i\ge1}$ is a sequence of independent and identically
distributed (i.i.d.) random elements taking values in $\mathbb{X}$ and
with probability
distribution $P_0$, the nonnegative random weights $(\tilde p_i)_{i\ge
1}$ are independent from $(Y_i)_{i\ge1}$ and are such that $\sum
_{i\ge
1} \tilde p_i=1$, almost surely. In the species sampling context, the
$Y_i$'s act as species tags and $\tilde p_i$ is the random proportion
with which the $i$th species is present in the population. If
$(X_n)_{n\ge1}$ is an exchangeable sequence directed by a species
sampling model $\tilde p$, that is, for every $n \ge1$ and $A_1,\ldots
,A_n$ in $\mathscr{X}$ one has
%
%
\begin{equation}
\label{eqexchange}
\mathbb{P}[X_1\in A_1,\ldots,X_n\in A_n | \tilde p]=\prod_{i=1}^n
\tilde p(A_i)
\end{equation}
almost surely, then $(X_n)_{n\ge1}$ is termed \textit{species sampling
sequence}.
Besides being an effective tool for statistical inference, species
sampling models
have an appealing structural property established in \cite{Pit96a}.
Indeed, if $(X_n)_{n\ge1}$ is a species sampling sequence, then there
exists a
collection of nonnegative weights $\{p_{j,n}(n_1,\ldots,n_k)\dvt 1\le
j\le k+1, 1\le k\le n, n\ge1\}$ such that $\sum
_{j=1}^{k+1}p_{j,n}(n_1,\ldots,n_k)=1$, for any vector of positive
integers $(n_1,\ldots,n_k)$ with $\sum_{j=1}^k n_j=n$, and
\[
\mathbb{P}[X_{n+1} \in \cdot| X_1,\ldots
,X_n]=p_{K_n+1,n}(n_1,\ldots
,n_{K_n})P_0( \cdot)+\sum_{j=1}^{K_n} p_{j,n}(n_1,\ldots,n_{K_n})
\delta_{X^*_j}( \cdot),
\]
where $X_1,\ldots,X_n$ is a sample with $K_n$ distinct values
$X_1^*,\ldots,X_{K_n}^*$. Statistical applications involving species
sampling models for different purposes than those of the present paper
are provided, for example, in \cite{navarra,navarra2,Lij07}.

The Bayesian nonparametric approach we undertake postulates that the
data are exchangeable and generated by a species sampling model. Then,
conditionally on the basic sample of size $n$, inference is to be made
on the number $K_m^{(n)}$ of new distinct species that will be observed
in the additional sample of size $m$. Interest lies in providing both a
point estimate and a measure of uncertainty, in the form of a credible
interval, for $K_m^{(n)}$ given $(X_1, \ldots, X_n)$. Since the
conditional distribution of $K_m^{(n)}$ becomes intractable for large
sizes $m$ of the additional sample, one is led to studying its limiting
behaviour as $m$ increases. Such asymptotic results, in addition to
providing useful approximations to the required estimators, are also of
independent theoretical interest since they provide useful insight on
the behaviour of the models we focus on. The only discrete random
probability measure for which a conditional asymptotic result, similar
to the one investigated in this paper, is known, is the two-parameter
Poisson--Dirichlet process, shortly denoted as $\operatorname{PD}(\sigma,\theta)$.
According to \cite{Pit96a}, a $\operatorname{PD}(\sigma,\theta)$ process is a
species sampling model characterized by\looseness=-1
%
%
\begin{equation}
\label{eqpesi0}
p_{K_n+1,i}(n_1,\ldots,n_{K_n})=\frac{\theta+K_n\sigma}{\theta+n} ,\qquad
p_{j,n}(n_1,\ldots,n_{K_n})=\frac{n_j-\sigma}{\theta+n}
\end{equation}\looseness=0
with $j=1,\ldots,K_n$, $\sigma\in(0,1)$ and $\theta>-\sigma$. In this
case, \cite{Fav09} provide a result describing the conditional
limiting behaviour of $K_m^{(n)}$. In the present paper, we focus on an
alternative species sampling model, termed \textit{normalized
generalized gamma process} in \cite{Lij07}. As we shall see in the
next section, it depends on two parameters $\sigma\in(0,1)$ and
$\beta
>0$ and, for the sake of brevity, is denoted by $\operatorname{NGG}(\sigma,\beta)$.
Moreover, it is characterized by
\begin{eqnarray}
\label{eqpesi}
p_{K_n+1,n}(n_1,\ldots,n_{K_n})
&=&\frac{\sigma}{n} \frac{\sum_{l=0}^{n}
{{n}\choose{l}}(-1)^l \beta^{l/\sigma}
\Gamma(K_n+1-{l}/{\sigma};
\beta)}{\sum_{l=0}^{n-1}
{{n-1}\choose{l}}(-1)^l \beta^{l/\sigma}
\Gamma(K_n-{l}/{\sigma};
\beta)},\\
\label{eqpesi1}
p_{j,n}(n_1,\ldots,n_{K_n})
&=&(n_j-\sigma)
\frac{\sum_{l=0}^{n}{{n}\choose{l}}(-1)^l \beta^{l/\sigma}
\Gamma(K_n-{l}/{\sigma}; \beta)}{n \sum_{l=0}^{n-1}
{{n-1}\choose{l}}(-1)^l \beta^{l/\sigma}
\Gamma(K_n-{l}/{\sigma}; \beta)}
\end{eqnarray}
for any $j\in\{1,\ldots,K_n\}$, where $\Gamma(a;x)$ is the incomplete
gamma function. The $\operatorname{NGG}(\sigma,\beta)$ process prior has gained some
attention in the Bayesian literature and it has proved to be useful for
various applications such as those considered, for example, in \cite
{Lij07,Arg1,Arg2,griffin1,griffin2,griffin3}. It is to be noted that
the $\operatorname{NGG}(\sigma,\theta)$ does not feature a posterior structure that is
as tractable as the one associated to the $\operatorname{PD}(\sigma, \theta)$ process
(see, e.g., \cite{Pit96a,Carlton,Lij08,jlp0}). Nonetheless, in
terms of practical implementation, it is possible to devise efficient
simulation algorithms that allow for a full Bayesian analysis within
models based on a $\operatorname{NGG}(\sigma,\beta)$ prior. See \cite{lp} for a
review of such algorithms.

In the present manuscript, we will specify the asymptotic
behaviour of $K_m^{(n)}$, given the basic sample, as $m$ diverges and
highlight the interplay between the conditional distributions of the
$\operatorname{PD}(\sigma,\theta)$ and the $\operatorname{NGG}(\sigma,\beta)$ processes. Since the
posterior characterization of a $\operatorname{NGG}(\sigma,\beta)$ process is far more
involved than the one associated to the $\operatorname{PD}(\sigma, \theta)$ process,
the derivation of the conditional asymptotic results considered in this
paper is technically more challenging. This is quite interesting since
it suggests that it is possible to study the limiting conditional
behaviour of $K_m^{(n)}$ even beyond species sampling models sharing
some sort of conjugacy property. For example, one might conjecture that
the same asymptotic regime, up to certain transformations of the
limiting random variable, should hold also for the wide class of
Gibbs-type priors, to be recalled in Section~\ref{section2}. An up to
date account of Bayesian Nonparametrics can be found in the monograph
\cite{hjort} and, in particular for asymptotic studies, \cite{ghosal}
provides a\vadjust{\goodbreak} review of asymptotics of nonparametric models in terms of
``frequentist consistency.'' Yet another type of asymptotic results are
obtained in \cite{pp,pecc}.

The outline of the paper is as follows. In Section~\ref{section2}, one can find a
basic introduction to species sampling models and a recollection of
some results in the literature concerning the asymptotic behaviour of
the number $K_{n}$ of distinct species in the basic sample, as $n$
increases. Section~\ref{section3} displays the main results, whereas the last
section contains some concluding remarks.

\section{Species sampling models and Gibbs-type priors}\label{section2}
Let us start by providing a succinct description of completely random
measures (CRM) before defining the specific models we will consider and
which can be derived as suitable transformations of CRMs. See \cite
{lp} for an overview of discrete nonparametric models defined in terms
of CRMs.

Suppose $\tilde\mu$ is a random element defined on some probability
space $(\Omega,\mathscr{F},\mathbb{P})$ and taking values on the
space $\mathcal{M}_\mathbb{X}
$ of boundedly finite measures on $(\mathbb{X}, \mathscr{X})$ such
that for any
$A_1,\ldots,A_n$ in $\mathscr{X}$, with $A_i\cap A_j=\varnothing$
for $i\ne
j$, the random
variables $\tilde\mu(A_1),\ldots,\tilde\mu(A_n)$ are mutually
independent. Then $\tilde\mu$ is termed \textit{completely random
measure} (CRM). It is well-known that the Laplace functional transform
of $\tilde\mu$ has a simple representation of the type
\[
\mathbb{E}[\mathrm{e}^{-\int f \,\mathrm{d}\tilde\mu}]=\mathrm
{e}^{-\psi(f)},
\]
where $\psi(f)=\int_{\mathbb{R}^+\times\mathbb{X}}[1-\mathrm
{e}^{-sf(y)}] \nu(\mathrm{d}s, \mathrm{d}
y)$ for any
measurable function $f\dvtx\mathbb{X}\to\mathbb{R}$ such that $\int|f|
\,\mathrm{d}\tilde\mu
<\infty$ almost surely and the measure $\nu$ on $\mathbb{R}^+\times
\mathbb{X}$ is known
as the L\'evy intensity of $\tilde\mu$.
See, for example, \cite{Kingman}. Since a CRM is almost surely
discrete, any CRM can be represented as
$\tilde\mu=\sum_{i\ge1} J_i \delta_{Y_i}$ with independent
random jump locations $(Y_i)_{i\ge1}$ and heights $(J_i)_{i\ge1}$.
For our purposes, it is enough to focus on the special case of $\nu$
factorizing as $\nu(\mathrm{d}s,\mathrm{d}x)=\rho(s) \,\mathrm{d}s
\alpha(\mathrm{d}x)$,
which implies independence of the locations $Y_i$'s and jumps $J_i$'s
in the above series representation. Furthermore, $\alpha$ can be taken
to be nonatomic and finite, the latter ensuring almost sure finiteness
of the corresponding CRM. Now, if $\operatorname{card}(\{J_i\dvtx  i\ge1\}\cap(0,\varepsilon
))=\int_0^\varepsilon\rho(s) \,\mathrm{d}s=\infty$ for any $\varepsilon
>0$, one can define a
random probability measure on $\mathbb{X}$ as
%
%
\begin{equation}
\label{eqnrmi}
\tilde p=\frac{\tilde\mu}{\tilde\mu(\mathbb{X})}.
\end{equation}
This family of random probability measures is known from \cite{jlp1} as
homogeneous normalized random measure with independent increments, a
subclass of the general class of normalized processes introduced in
\cite{rlp}. Note that an $\mathbb{X}$-valued exchangeable sequence
$(X_n)_{n\ge1}$ generated by $\tilde p$ as in \eqref{eqnrmi} is a
species sampling sequence.

Here we focus on a specific example where the CRM defining $\tilde p$
in \eqref{eqnrmi}
is the so-called \textit{generalized gamma process} \cite{brix} that
is characterized by
\[
\label{eqlevygg}
\rho(s)=\frac{\sigma}{\Gamma(1-\sigma)} s^{-1-\sigma} \mathrm
{e}^{-\tau s}\vadjust{\goodbreak}
\]
with $\sigma\in(0,1)$ and $\tau>0$. In this case,
%
%
\begin{equation}
\psi(f)=\int_\mathbb{X}\bigl[\bigl(f(x)+\tau\bigr)^\sigma-\tau^\sigma\bigr] \alpha
(\mathrm{d}x)
\label{eqlaplexpon}
\end{equation}
for any measurable function $f\dvtx \mathbb{X}\to\mathbb{R}$ such that
$\int|f|^\sigma\,\mathrm{d}
\alpha<\infty$. In the sequel the model will be reparameterized,
without loss of generality (see, e.g., \cite{Pit03,Lij07}), by
setting $\beta:=\tau^\sigma$ and $\alpha$ as a probability measure. The
corresponding CRM will be denoted by $\tilde\mu_{\sigma,\beta}$.
Henceforth, the random probability measure $\tilde p$ obtained by
normalizing $\tilde\mu_{\sigma,\beta}$ as in \eqref{eqnrmi} coincides,
in distribution, with the $\operatorname{NGG}(\sigma, \beta)$ process prior. An
important special case arises when $\beta=0$, since $\tilde{\mu
}_{\sigma
,0}$ reduces to the $\sigma$-stable process, which plays a key role
within the paper. For example, it is worth noting that $\tilde\mu
_{\sigma,\beta}$ can also be defined as
an exponential tilting of $\tilde\mu_{\sigma,0}$, for any $\beta>0$.
Specifically, if $\mathbb{P}_{\sigma,0}$ is
the probability distribution of $\tilde\mu_{\sigma,0}$ on $\mathcal
{M}_\mathbb{X}$ and $\mathbb{P}_{\sigma,\beta}$
is a probability measure on $\mathcal{M}_\mathbb{X}$ that is absolutely
continuous with
respect to $\mathbb{P}_{\sigma,0}$ and such that
%
%
\begin{equation}
\frac{\mathrm{d}\mathbb{P}_{\sigma,\beta}}{\mathrm{d}\mathbb
{P}_{\sigma,0}}(\mu)=
\exp\{\beta-\beta^{1/\sigma} \mu(\mathbb{X})\}
\label{eqexpontilt}
\end{equation}
then $\mathbb{P}_{\sigma,\beta}$ coincides with the probability
distribution of
$\tilde\mu_{\sigma,\beta}$. In
a similar fashion, one can also define the $\operatorname{PD}(\sigma,\theta)$ process
as a polynomial tilting
of $\tilde\mu_{\sigma,0}$, for any $\theta>-\sigma$. Indeed, one
introduces another probability measure $\mathbb{Q}_{\sigma,\theta}$
that is still absolutely continuous with respect to $\mathbb
{P}_{\sigma,0}$ and
whose Radon--Nykodim derivative is
%
%
\begin{equation}
\label{eqpolytilt}
\frac{\mathrm{d}\mathbb{Q}_{\sigma,\theta}}{\mathrm{d}\mathbb
{P}_{\sigma,0}}(\mu)=
\frac{\Gamma(\theta+1)}
{\Gamma({\theta}/{\sigma}+1)} [\mu(\mathbb{X})]^{-\theta}
\end{equation}
for any $\sigma\in(0,1)$ and $\theta>-\sigma$. If $\mu_{\sigma
,\theta
}^*$ is the random measure with probability distribution $\mathbb
{Q}_{\sigma,\theta}$ above, then $p^*=\mu_{\sigma,\theta}^*/\mu
_{\sigma
,\theta}^*(\mathbb{X})$ coincides,
in distribution, with a $\operatorname{PD}(\sigma,\theta)$ process. See \cite
{Pit06}. The different tilting structure
featured by the normalized generalized gamma process and the two parameter
Poisson--Dirichlet process will be reflected by the limiting results to
be illustrated
in the paper.

It is also worth to recall that both the $\operatorname{NGG}(\sigma,\beta)$ and the
$\operatorname{PD}(\sigma,\theta)$ processes can be seen as elements of the
general class of \textit{Gibbs-type} nonparametric priors introduced in~\cite{gnedin}.
Gibbs-type priors represent the most tractable subclass
of species sampling models. They are characterized by a parameter
$\sigma<1$ and a collection of non-negative quantities $\{V_{n,k}\dvt
n\ge1, 1\le k\le n\}$ that satisfy the forward recursive relations
\[
\label{eqgibbsrecurs}
V_{n,k}=V_{n+1,k+1}+(n-k\sigma)V_{n+1,k}.
\]
These $V_{n,k}$'s define the predictive weights that characterize a
species sampling
sequence governed by a Gibbs-type prior. Indeed, one has
%
%
\begin{equation}
\label{eqweightsgibbs}
p_{K_n+1,n}(n_1,\ldots,n_{K_n})=\frac{V_{n+1,K_n+1}}{V_{n,K_n}},\qquad
p_{j,n}(n_1,\ldots,n_{K_n})=\frac{V_{n+1,K_n}}{V_{n,K_n}}(n_j-\sigma)
\end{equation}
for each $j\in\{1,\ldots,K_n\}$. The fundamental simplification
involved in \eqref{eqweightsgibbs} is that the probability of
observing a ``new'' or an ``old'' species depend on the sample size and
on the number of already observed distinct species but not on their
frequencies: this crucially simplifies explicit calculations. Turning
to the two specific processes introduced before, in accordance with
\eqref{eqpesi0}, the $\operatorname{PD}(\sigma,\theta)$ process identifies a
Gibbs-type prior with
\[
\label{eqvnkpd}
V_{n,k}=\frac{\prod_{i=1}^{k-1}(\theta+i\sigma)}{(\theta+1)_{n-1}},
\]
whereas, in accordance with \eqref{eqpesi} and \eqref{eqpesi1}, a
$\operatorname{NGG}(\sigma,\beta)$ prior is also of Gibbs-type with $\sigma\in
(0,1)$ and
\[
\label{eqvnkngg}
V_{n,k}=\frac{\mathrm{e}^\beta\sigma^{k-1}}{\Gamma(n)} \sum_{l=0}^{n-1}
\pmatrix{{n-1}\cr{l}}(-1)^{l} \beta^{l/\sigma} \Gamma\biggl(k-\frac{l}{\sigma
};\beta\biggr).
\]
As shown in \cite{Lij108}, a normalized CRM is a Gibbs-type prior (with
$\sigma\in(0,1)$) if and only if it is a $\operatorname{NGG}(\sigma,\beta)$ process.
This result also motivates
the focus of the paper on the $\operatorname{NGG}(\sigma, \beta)$ process, which
clearly has a prominent role.

The result on the limiting behaviour of $K_m^{(n)}$ to be determined in
the next section parallels known results for the unconditional case
where one aims at determining the asymptotics of $K_n$ as the sample
size $n$ increases and connects to the conditional asymptotics
displayed in \cite{Fav09} for the $\operatorname{PD}(\sigma,\theta)$ process. In
order to describe the result for the unconditional case, let $T_{\sigma
,0}:=\tilde\mu_{\sigma,0}(\mathbb{X})$ be the random total mass of
a $\sigma$-stable CRM and denote by $f_\sigma$ its density function which
satisfies $\int_0^\infty\mathrm{e}^{-\lambda s} f_\sigma(s) \,\mathrm
{d}s=\mathrm{e}
^{-\lambda^\sigma}$
for any $\lambda>0$. Moreover, let $T_{\sigma,\beta}:=\tilde\mu
_{\sigma
,\beta}(\mathbb{X})$ be the random total mass of $\operatorname{NGG}(\sigma, \beta
)$ process
and recall that its law can be obtained by exponentially
tilting the probability distribution of $T_{\sigma,0}$ as in \eqref
{eqexpontilt}. In particular,
if
%
%
\begin{equation}\label{prima}
S_{\sigma,\beta}\stackrel{\mathrm{d}}{=}T_{\sigma,\beta}^{-\sigma},
\end{equation}
then its density function, with respect to the Lebesgue measure on
$\mathbb{R}$, coincides
with
\[
\label{eqdensdiversity}
g_{\sigma,\beta}(s)=
\frac{\mathrm{e}^\beta}{\sigma} \mathrm{e}^{-({\beta
}/{s})^{
{1}/{\sigma}}}
s^{-1-{1}/{\sigma}} f_\sigma(s^{-1/\sigma}) \mathbh
{1}_{(0,\infty)}(s)
\]
and one has that
%
%
\begin{equation}
\label{eqdiversity}
\frac{K_n}{n^\sigma} \stackrel{\mathrm{a.s.}}{\longrightarrow}
S_{\sigma
,\beta}.
\end{equation}
According to the terminology introduced by \cite{Pit03}, the random variable
$S_{\sigma,\beta}$ is the so-called $\sigma$-\textit{diversity} of
the exchangeable random
partition induced by a $\operatorname{NGG}(\sigma,\beta)$ process prior. See also
Definition 3.10 in
Pitman \cite{Pit06}. Note that a similar result holds true for the
$\operatorname{PD}(\sigma,\theta)$ process.
Indeed, if $T^{\prime}_{\sigma,\theta}\stackrel{\mathrm{d}}{=}\mu
^*_{\sigma
,\theta}(\mathbb{X})$ so that its probability
distribution is obtained by polynomially tilting the probability
distribution of
$T_{\sigma,0}$ as in \eqref{eqpolytilt} and
%
%
\begin{equation}\label{seconda}
S_{\sigma,\theta}'\stackrel{\mathrm{d}}{=}(T_{\sigma,\theta
}')^{-\sigma}
\end{equation}
admits density function
\[
\label{eqdensdiversitypd}
h_{\sigma,\theta}(s)=\frac{\Gamma(\theta+1)}{\Gamma({\theta
}/{\sigma}+1)}
\frac{s^{{\theta}/{\sigma}-{1}/{\sigma}-1}}{\sigma}
f_\sigma
(s^{-1/\sigma}) \mathbh{1}_{(0,\infty)}(s).
\]
Then one has
%
%
\begin{equation}
\label{eqdiversitypd}
\frac{K_n}{n^\sigma} \stackrel{\mathrm{a.s.}}{\longrightarrow}
S_{\sigma
,\theta}'.
\end{equation}
See \cite{Pit06}, Theorem 3.8. These results are somehow
in line with the fact that the $\operatorname{NGG}(\sigma,\beta)$ and the $\operatorname{PD}(\sigma
,\theta)$ processes are distributionally
equivalent to normalized random measures that are obtained by an
exponential and
a polynomial tilting, respectively, of a $\sigma$-stable CRM as
highlighted in
\eqref{eqexpontilt} and in \eqref{eqpolytilt}. Finally, note that a
combination of \cite{gnedin}, Theorem 12, and \cite{Pit03}, Proposition
13,
shows that the unconditional asymptotic results in \eqref
{eqdiversity} and
\eqref{eqdiversitypd} can be extended to the whole class of
Gibbs-type priors. See also \cite{griffiths07} for another
contribution at the interface between Bayesian Nonparametrics and
Gibbs-type random partitions.

\section{\texorpdfstring{Asymptotics of $K_m^{(n)}$ with a $\operatorname{NGG}(\sigma,\beta)$ process}
{Asymptotics of K m (n) with a NGG(sigma,beta) process}}\label{section3}

As mentioned before, inference on $K_m^{(n)}$ is of great importance
since it provides a measure of species richness of a community of
plants/animals or of a cDNA library for gene discovery. The key
quantity for obtaining posterior inferences is given by the probability
distribution $\mathbb{P}[K_{m}^{(n)}=k | X_1, \ldots, X_n]$ for $k=0,
\ldots, m$. By virtue of predictive sufficiency of the number $K_n$ of
distinct species observed among the first $n$ data $X_1,\ldots,X_n$, in
\cite{Lij207} it has been shown that in the $\operatorname{NGG}(\sigma, \beta)$
this distribution coincides with
\begin{eqnarray}
\label{posteriordistinct}
P_{m}^{(n,j)}(k)
&:=&
\mathbb{P}\bigl[K_{m}^{(n)}=k | K_{n}=j\bigr]
\nonumber
\\[-8pt]
\\[-8pt]
\nonumber
&=&
\frac{\mathscr{G}(m,k;\sigma,-n+j\sigma)}{(n)_m}
\frac{\sum_{l=0}^{n+m-1}{n+m-1\choose l}(-1)^{l}\beta^{l/\sigma
}\Gamma
(j+k-{l}/{\sigma};\beta)}{\sum_{l=0}^{n-1}{n-1\choose
l}(-1)^{l}\beta^{l/\sigma}\Gamma(j-{l}/{\sigma};\beta)}
\end{eqnarray}
for $k=0,\ldots,m$, with $\mathscr{G}(n,k;s,r)$ denoting the
non--central generalized factorial coefficient. See \cite{Cha05} for
a comprehensive account on generalized factorial coefficients.
Expression \eqref{posteriordistinct} can be interpreted as the
``posterior'' probability distribution of the number of distinct new
species to be observed in a further sample of size $m$.
Now, based on \eqref{posteriordistinct}, one obtains the expected
number of new species as
%
%
\begin{equation}
\label{posteriormean}
\hat{E}_{m}^{(n,j)}:=\mathbb{E}\bigl[K_{m}^{(n)}|K_{n}=j\bigr]=\sum
_{k=0}^{m}kP_{m}^{(n,j)}(k),
\end{equation}
which corresponds to the Bayes estimator of $K_{m}^{(n)}$ under
quadratic loss. Moreover, a measure of uncertainty of the point
estimate $\hat{E}_{m}^{(n,j)}$ can be obtained in terms of $\alpha
$-credible intervals that is, by determining an interval $(z_1, z_2)$
with $z_1<z_2$ such that $\mathbb{P}[z_1\leq K_{m}^{(n)}\leq z_2 | K_{n}=j]
\geq\alpha$. The interval $(z_1, z_2)$ of shortest length is then
typically referred to as highest posterior density interval.

The main advantage of the distribution \eqref{posteriordistinct} is
that it is explicit. However, since the sum of incomplete gamma
functions cannot be further simplified, its computation can become
overwhelming even for moderately large sizes of $n$ and $m$. This fact
represents a major problem in the frequent practical situations in
which the size of the additional sample of interest is large. For
instance, in genomic applications one has to deal with relevant
portions of cDNA libraries which typically consist of millions of
genes. Hence, it is natural to study the asymptotics for $K_{m}^{(n)}$,
given $K_{n}$, as $m\rightarrow+\infty$, in order to obtain
approximations of \eqref{posteriordistinct} and, consequently, also of
\eqref{posteriormean} and of the corresponding highest posterior
density intervals. Indeed, if one is able to show that a suitable
rescaling of $K_{m}^{(n)}$, given $K_{n}$, converges in law to some
random variable, one can use the probability distribution of this
limiting random quantity in order to derive the desired approximations.

\subsection{Asymptotic distribution}
The statement of the main result in the paper involves a positive
random variable $Y_q$ whose density function is, for any $q>0$,
\[
\label{eqmittag}
f_{Y_q}(y)=\frac{\Gamma(q\sigma+1)}{\sigma\Gamma(q+1)}
y^{q-1-1/\sigma} f_\sigma(y^{-1/\sigma})
\]
and we $B_{a,b}$ to denote a beta random variable with parameters
$(a,b)$. Moreover, set $S_{n,j}\stackrel{\mathrm{d}}{=} B_{j,n/\sigma
-j}Y_{n/\sigma}$, with $B_{j,n/\sigma-j}$ and $Y_{n/\sigma}$
independent, and denote by $g_{S_{n,j}}$ the density function of $S_{n,j}$.

\begin{thm}\label{limitresult}
If $(X_n)_{n\ge1}$ is a species sampling sequence directed by a
$\operatorname{NGG}(\sigma,\beta)$ process prior, conditional on $K_{n}=j$ one has
%
%
\begin{equation}\label{result}
\frac{K_{m}^{(n)}}{m^{\sigma}}\rightarrow Z_{n,j}\qquad \mbox{a.s.}
\end{equation}
as $m\rightarrow+\infty$, where $Z_{n,j}$ is a positive random variable
obtained by exponentially tilting the density function of $S_{n,j}$, namely
\[
f_{Z_{n,j}}(z)=\frac{\Gamma(j)\mathrm{e}^{-(\beta/z)^{1/\sigma
}}g_{S_{n,j}}(z)}{\sum_{l=0}^{n-1}{n-1\choose l}(-1)^{l}\beta
^{l/\sigma
}\Gamma(j-{l}/{\sigma};\beta)}.
\]
\end{thm}

\begin{pf}
The first part of the proof exploits a martingale convergence theorem
along the same lines of \cite{Pit06}, Theorem 3.8. In particular,
let us start by computing the likelihood ratio
\[
\label{likeratio}
M^{(n)}_{\sigma,\beta,m}=\frac{\mathrm{d}\mathbb{P}^{(n)}_{\sigma
,\beta}}{\mathrm{d}\mathbb{P}
^{(n)}_{\sigma,0}}\bigg|_{\mathscr{F}^{(n)}_{m}}=\frac{q^{(n)}_{\sigma
,\beta}(K_{m}^{(n)})}{q_{\sigma,0}^{(n)}(K_{m}^{(n)})},
\]
where $\mathscr{F}^{(n)}_{m}=\sigma(X_{n+1},\ldots,X_{n+m})$,
$\mathbb{P}
_{\sigma,\beta}^{(n)}$ is the conditional
probability distribution of a normalized generalized gamma process with
parameter $(\sigma,\beta)$ given
$K_{n}$ and, by virtue of \cite{Lij08}, Proposition 1,
\[
q^{(n)}_{\sigma,\beta}(K_{m}^{(n)})
=\frac{\sigma^{K_{m}^{(n)}}}{(n)_m} \frac{\sum
_{l=0}^{n+m-1}{n+m-1\choose l}(-1)^{l}\beta^{l/\sigma}\Gamma
(K_{n}+K_{m}^{(n)}-{l}/{\sigma};\beta)}
{\sum_{l=0}^{n-1}{n-1\choose l}(-1)^{l}\beta^{l/\sigma}\Gamma
(K_{n}-{l}/{\sigma};\beta)}
\]
for any integer $K_{n}\geq1$ and
\[
q^{(n)}_{\sigma,0}(K_{m}^{(n)})=\frac{\sigma^{K_{m}^{(n)}}
(K_n)_{K_m^{(n)}}}{(n)_m}.
\]
Hence, $(M^{(n)}_{\sigma,\beta,m},\mathscr{F}^{(n)}_{m})_{m\geq1 }$ is
a $\mathbb{P}^{(n)}_{\sigma,0}$-martingale and by a martingale
convergence theorem, $M^{(n)}_{\sigma,\beta,m}$ has a $\mathbb
{P}^{(n)}_{\sigma
,0}$ almost sure limit, say $M^{(n)}_{\sigma,\beta}$, as
$m\rightarrow+\infty$. Clearly, we have that \mbox{$\mathbb{E}_{\sigma
,0}^{(n)}[M^{(n)}_{\sigma,\beta}]=1$}, where $\mathbb{E}_{\sigma,0}^{(n)}$
denotes the expected
value with respect to $\mathbb{P}_{\sigma,0}^{(n)}$. Let now
$(E_{n})_{n\geq1}$
be a sequence of i.i.d. random variables having a negative exponential
distribution with parameter~$1$. Moreover,
suppose the $E_{n}$'s are independent of $(K_{n},K_{m}^{(n)})$. Set
$\mathscr{E}_{m}^{(n)}:=\sum_{i=1}^{K_{n}+K_{m}^{(n)}}E_{i}$
and note that, conditionally on $(K_{n},K_{m}^{(n)})$, $\mathscr
{E}^{(n)}_{m}$ has gamma distribution
with expected value $K_{n}+K_{m}^{(n)}$. We can then rewrite
$M^{(n)}_{\sigma,\beta,m}$ as follows
\begin{eqnarray*}
M^{(n)}_{\sigma,\beta,m}
&=&\frac{\Gamma(K_{n})}{\sum_{l=0}^{n-1}{n-1\choose l}(-1)^{l}\beta
^{l/\sigma}\Gamma(K_{n}-{l}/{\sigma};\beta)}
\frac{1}{\Gamma(K_{n}+K_{m}^{(n)})}\\
&&{} \times\sum_{l=0}^{n+m-1}{n+m-1\choose
l}(-1)^{l}\beta^{l/\sigma}\int_{\beta}^{+\infty
}y^{K_{n}+K_{m}^{(n)}-l/\sigma-1}\mathrm{e}^{-y}\,\mathrm{d}y\\
&=&\frac{\Gamma(K_{n})}{\sum_{l=0}^{n-1}{n-1\choose l}(-1)^{l}\beta
^{l/\sigma}\Gamma(K_{n}-{l}/{\sigma};\beta)}
\frac{1}{\Gamma(K_{n}+K_{m}^{(n)})}\\
&&{} \times\int_{\beta}^{+\infty
}y^{K_{n}+K_{m}^{(n)}-1}\mathrm{e}^{-y}\biggl(1-\frac{\beta^{1/\sigma
}}{y^{1/\sigma}}\biggr)^{n+m-1}\,\mathrm{d}y\\
&=&\frac{\Gamma(K_{n})}{\sum_{l=0}^{n-1}{n-1\choose l}(-1)^{l}\beta
^{l/\sigma}\Gamma(K_{n}-{l}/{\sigma};\beta)}\\
&&{}
\times\mathbb{E}\biggl[\mathbh{1}_{(\beta,+\infty)}\bigl(\mathscr
{E}_{m}^{(n)}\bigr)\biggl(1-\frac{\beta^{1/\sigma}}{(\mathscr
{E}_{m}^{(n)})^{1/\sigma}}\biggr)^{n+m+1}\bigg|\mathscr
{F}_{m}^{(n)}\biggr].
\end{eqnarray*}
From the strong law of large numbers, $\mathscr
{E}^{(n)}_{m}/(K_{n}+K_{m}^{(n)})\rightarrow1$ as $m\rightarrow
+\infty$
and conditionally on $(K_{n},K_{m}^{(n)})$. Using the dominated
convergence theorem, we have
\begin{eqnarray*}
M^{(n)}_{\sigma,\beta,m}&\approx&\frac{\Gamma(K_{n})}{\sum
_{l=0}^{n-1}{n-1\choose l}(-1)^{l}\beta^{l/\sigma}\Gamma
(K_{n}-{l}/{\sigma};\beta)}\\
&&{}
\times\biggl(1-\frac{\beta^{1/\sigma} }{( (K_{n}+K_{m}^{n})(\mathscr
{E}^{n}_{m}/(K_{n}+K_{m}^{n})))^{1/\sigma}}\biggr)^{n+m-1}\\[-2pt]
&\approx&\frac{\Gamma(K_{n})}{\sum_{l=0}^{n-1}{n-1\choose
l}(-1)^{l}\beta^{l/\sigma}\Gamma(K_{n}-{l}/{\sigma};\beta
)}\biggl(1-\frac{\beta^{1/\sigma}}{(K_{n}+K_{m}^{(n)})^{1/\sigma}}
\biggr)^{n+m-1}\\[-2pt]
&\approx&\frac{\Gamma(K_{n})}{\sum_{l=0}^{n-1}{n-1\choose
l}(-1)^{l}\beta^{l/\sigma}\Gamma(K_{n}-{l}/{\sigma};\beta
)}
\exp\biggl\{-m\frac{\beta^{1/\sigma}}{(K_{m}^{(n)})^{1/\sigma}}\biggr\}
\end{eqnarray*}
as $m\rightarrow+\infty$. Since $M^{(n)}_{\sigma,\beta
,m}\rightarrow
M^{(n)}_{\sigma,\beta}$ almost surely (with respect to $\mathbb
{P}^{(n)}_{\sigma
,0}$), then
there exists some positive random variable, say $L_{\sigma,n}$ such
that $m/ (K_{m}^{(n)})^{1/\sigma}\rightarrow L_{\sigma,n}$
almost surely (with respect to $\mathbb{P}^{(n)}_{\sigma,0}$). In
order to
identify the probability distribution of $L_{\sigma,n}$, note that it
must be such that
%
%
\begin{equation}\label{condition}
\mathbb{E}[\mathrm{e}^{-\beta^{1/\sigma}L_{\sigma,n}}]
=
\frac{1}{\Gamma(K_{n})}\int_{\beta}^{+\infty}y^{K_{n}-1}\biggl(1-\frac
{\beta^{1/\sigma}}{y^{1/\sigma}}
\biggr)^{n-1}\mathrm{e}^{-y} \,\mathrm{d}y.
\end{equation}
Since $S_{n,K_n}\stackrel{\mathrm{d}}{=}B_{K_{n},n/\sigma
-K_{n}}Y_{n/\sigma}$, we have to prove that $L_{\sigma,n}\stackrel
{\mathrm
{d}}{=}S_{n,K_n}^{-1/\sigma}$, that is, that the density function of
$L_{\sigma,n}$ coincides with
%
\begin{eqnarray}
\label{density1}
f_{L_{\sigma,n}}(z)
&=&\frac{\sigma\Gamma(n)}{\Gamma(K_{n})\Gamma({n}/{\sigma
}-K_{n})} z^{-\sigma-1}
\nonumber
\\[-9pt]
\\[-9pt]
\nonumber
&&{}\times\int_{z^{-\sigma}}^{+\infty}\frac{1}{v}v^{{n}/{\sigma
}-1-{1}/{\sigma}}f_{\sigma}(v^{-{1}/{\sigma}})
\biggl(\frac{z^{-\sigma}}{v}\biggr)^{K_{n}-1}\biggl(1-\frac{z^{-\sigma
}}{v}\biggr)^{{n}/{\sigma}-K_{n}-1} \,\mathrm{d}v.\qquad
\end{eqnarray}
So we simply have to show that the Laplace transform of the density
function in \eqref{density1} is given by~\eqref{condition}. By a
simple change of variable, $x=v^{-1/\sigma}$, the previous density
reduces~to
\[
f_{L_{\sigma,n}}(z)=\frac{\sigma^{2}\Gamma(n)}{\Gamma(K_{n})\Gamma
({n}/{\sigma}-K_{n})}\int_{0}^{z}x^{-n+\sigma}f_{\sigma}
(x)\biggl(\frac{z^{-\sigma}}{x^{-\sigma}}\biggr)^{K_{n}-1}
\biggl(1-\frac{z^{-\sigma}}{x^{-\sigma}}\biggr)^{{n}/{\sigma
}-K_{n}-1}z^{-\sigma-1}\,\mathrm{d}x
\]
and by the change of variable $y=z^{-\sigma}/x^{-\sigma}$
\[
f_{L_{\sigma,n}}(z)=\frac{\sigma\Gamma(n)}{\Gamma(K_{n})\Gamma
({n}/{\sigma}-K_{n})}z^{-n} \int_{0}^{1}y^{-n/\sigma+1/\sigma
+K_{n}-1}(1-y)^{n/\sigma-K_{n}-1}f_{\sigma}(zy^{1/\sigma}
)\,\mathrm{d}y.
\]
The Laplace transform of $f_{L_{\sigma,n}}$ is then given by
\begin{eqnarray*}
\mathbb{E}[\mathrm{e}^{-\beta^{1/\sigma}L_{\sigma,n}}]
&=&\frac{\sigma\Gamma(n)}{\Gamma(K_{n})\Gamma(
{n}/{\sigma
}-K_{n})}\\[-2pt]
&&{}\times \int_{0}^{\infty}\mathrm{e}^{\beta^{1/\sigma}z}z^{-n}\\[-2pt]
&&\hspace*{28pt}{} \times\int_{0}^{1}y^{-n/\sigma
+1/\sigma+K_{n}-1}(1-y)^{n/\sigma-K_{n}-1}f_{\sigma}(zy^{1/\sigma
}) \,\mathrm{d}y \,\mathrm{d}z \\[-2pt]
&=&\frac{\sigma\Gamma(n)}{\Gamma(K_{n})\Gamma(
{n}/{\sigma
}-K_{n})}\int_{0}^{1}y^{-n/\sigma+1/\sigma+K_{n}-1}(1-y)^{n/\sigma
-K_{n}-1}\\[-2pt]
&&\hspace*{99pt}{} \times\int_{0}^{\infty}
\mathrm{e}^{\beta^{1/\sigma}z}z^{-n}f_{\sigma}(zy^{1/\sigma}
)\,\mathrm{d}z\,\mathrm{d}y\\[-2pt]
&=&\frac{\sigma\Gamma(n)}{\Gamma(K_{n})\Gamma(
{n}/{\sigma
}-K_{n})}
\int_{0}^{1}y^{K_{n}-1}(1-y)^{n/\sigma-K_{n}-1}\\[-2pt]
&&\hspace*{99pt}{} \times
\int_{0}^{\infty}\mathrm{e}^{(\beta/y)^{1/\sigma}h}h^{-n}f_{\sigma}
(h)\,\mathrm{d}h \,\mathrm{d}y\\[-2pt]
&=&\frac{\Gamma({n}/{\sigma})}{\Gamma(K_{n})\Gamma
({n}/{\sigma}-K_{n})}\int
_{0}^{1}y^{K_{n}-1}(1-y)^{n/\sigma-K_{n}-1}\\[-2pt]
& &\hspace*{99pt}{}\times
\frac{\sigma\Gamma(n)}{\Gamma({n}/{\sigma})
}\int_{0}^{\infty}\mathrm{e}^{(\beta/y)^{1/\sigma
}h}h^{-n}f_{\sigma}
(h)\,\mathrm{d}h
\,\mathrm{d}y.
\end{eqnarray*}
According to the well-known gamma identity, we can write
\[
\frac{\sigma\Gamma(n)}{\Gamma({n}/{\sigma})}\int
_{0}^{\infty}
\frac{\mathrm{e}^{(\beta/y)^{1/\sigma}h}}{h^{n}} f_{\sigma}(h
) \,\mathrm{d}h
=\frac{\sigma}{\Gamma({n}/{\sigma})}
\int_{0}^{\infty}u^{n-1}\int_{0}^{\infty}\mathrm{e}^{-h(
{\beta^{1/\sigma
}}/{y^{1/\sigma}}+u)}
f_{\sigma}(h) \,\mathrm{d}h \,\mathrm{d}u
\]
obtaining
\begin{eqnarray*}
&&\frac{\Gamma({n}/{\sigma})}{\Gamma(K_{n})\Gamma
({n}/{\sigma}-K_{n})}\int_{0}^{1}y^{K_{n}-1}(1-y)^{n/\sigma
-K_{n}-1}\\[-2pt]
& &\qquad{}\times
\frac{\sigma}{\Gamma({n}/{\sigma})}
\int_{0}^{+\infty}u^{n-1}\int_{0}^{+\infty}\mathrm{e}^{-h(
{\beta
^{1/\sigma}}/{y^{1/\sigma}}+u)} f_{\sigma}(h)\,\mathrm{d}h \,\mathrm
{d}u\\[-2pt]
&&\quad=\frac{\Gamma({n}/{\sigma})}{\Gamma(K_{n})\Gamma
({n}/{\sigma}-K_{n})}\int_{0}^{1}y^{K_{n}-1}(1-y)^{n/\sigma
-K_{n}-1}\\[-2pt]
& &\qquad\times
\frac{\sigma}{\Gamma({n}/{\sigma})}\int_{0}^{+\infty
}u^{n-1} \mathrm{e}^{-({\beta^{1/\sigma}}/{y^{1/\sigma
}}+u)^{\sigma}}
\,\mathrm{d}u\\[-2pt]
&&\quad=
\frac{\Gamma({n}/{\sigma})}{\Gamma(K_{n})\Gamma
({n}/{\sigma}-K_{n})}\int_{0}^{1}y^{K_{n}-1}(1-y)^{n/\sigma
-K_{n}-1}\\[-2pt]
&&\qquad{} \times\frac{1}{\Gamma(
{n}/{\sigma})}\int_{\beta}^{+\infty}z^{n/\sigma-1}\biggl(1-
\biggl(\frac{\beta}{zy}\biggr)^{1/\sigma}\biggr)^{n-1}\mathrm{e}^{-z} \,\mathrm{d}z
\,\mathrm{d}y\\[-2pt]
&&\quad=\frac{1}{\Gamma(K_{n})\Gamma({n}/{\sigma}-K_{n})}\\[-2pt]
&&\qquad{}\times\int
_{\beta}^{+\infty}\mathrm{e}^{-z} \\[-2pt]
& &\qquad{}\hspace*{36pt}\times\int
_{0}^{z}w^{K_{n}-1}(z-w)^{n/\sigma-K_{n}-1}\biggl(1-\biggl(\frac{\beta
}{w}\biggr)^{1/\sigma}\biggr)^{n-1} \,\mathrm{d}w \,\mathrm{d}z\\[-2pt]
&&\quad=\frac{1}{\Gamma(K_{n})\Gamma({n}/{\sigma}-K_{n})}\int
_{\beta}^{\infty}w^{K_{n}-1}\biggl(1-\biggl(\frac{\beta}{w}
\biggr)^{1/\sigma}\biggr)^{n-1}\\[-2pt]
&&\hspace*{102pt}\qquad{} \times\int_{w}^{+\infty}\,\mathrm
{e}^{-z}(z-w)^{n/\sigma-K_{n}-1}
\,\mathrm{d}z \,\mathrm{d}w
\end{eqnarray*}
which corresponds to \eqref{condition}. Finally, since the probability
measures $\mathbb{P}_{\beta,\sigma}^{(n)}$ and
$\mathbb{P}_{0,\sigma}^{(n)}$ are mutually absolutely continuous,
almost sure
convergence holds true with respect
to $\mathbb{P}_{\beta,\sigma}^{(n)}$, as well. In order to deduce
the $\mathbb{P}_{\beta
,\sigma}^{(n)}$-law of $Z_{n,K_{n}}$, it is sufficient
to exploit a change of measure suggested by
\[
\mathbb{P}^{(n)}_{\sigma,\beta}(A)=\int_{A}\frac{\mathrm
{d}\mathbb{P}_{\sigma,\beta
}^{(n)}}{\mathrm{d}\mathbb{P}_{\sigma,0}^{(n)}}\,\mathrm{d}
\mathbb{P}_{\sigma,0}^{(n)}
\]
and by the fact that
\[
\frac{\mathrm{d}\mathbb{P}_{\sigma,\beta}^{(n)}}{\mathrm
{d}\mathbb{P}_{\sigma
,0}^{(n)}}=M^{(n)}_{\sigma,\beta}=\frac{\Gamma(K_{n})}{\sum
_{l=0}^{n-1}{n-1\choose l}(-1)^{l}\beta^{l/\sigma}\Gamma
(K_{n}-{l}/{\sigma};\beta)}\mathrm{e}^{-\beta^{1/\sigma
}L_{\sigma,n}}.
\]
This completes the proof.
\end{pf}

It is worth stressing that the limit random variable in the conditional
case is the same as in the unconditional case but with updated
parameters and a rescaling induced by a beta-distributed random
variable. The density of $Z_{n,j}$ in \eqref{result} can be formally
represented as
\begin{eqnarray}\label{density}
f_{Z_{n,j}}(z)&=&\frac{\Gamma(j)\mathrm{e}^{-(\beta/z)^{1/\sigma
}}}{\sum
_{l=0}^{n-1}{n-1\choose l}(-1)^{l}\beta^{l/\sigma}\Gamma(j-
{l}/{\sigma};\beta)}
\nonumber
\\[-8pt]
\\[-8pt]
\nonumber
&&{}
\times\frac{\Gamma(n)}{\Gamma(j)\Gamma({n}/{\sigma}-j
)}z^{j-1}\int_{z}^{+\infty}v^{-1/\sigma}(v-z)^{n/\sigma
-j-1}f_{\sigma
}(v^{-1/\sigma})\,\mathrm{d}v,
\end{eqnarray}
where we recall that $f_\sigma$ is the density function of a positive
stable random variable and, then, coincides with a density function of
the random total mass of a $\sigma$-stable CRM $T_{\sigma,0}:=\tilde
\mu
_{\sigma,0}(\mathbb{X})$. Theorem~\ref{limitresult} can be compared
with an
analogous result recently obtained in~\cite{Fav09}, Proposition 2,
for the $\operatorname{PD} (\sigma,\theta)$ process, where it is shown that the number
of new distinct species $K_{m}^{(n)}$ induced by the $\operatorname{PD} (\sigma,\theta
)$ process is such that
%
%
\begin{equation}\label{limittwoparameter}
\frac{K_{m}^{(n)}}{m^{\sigma}}\stackrel{\mathrm
{a.s.}}{\longrightarrow}
Z^{\prime}_{n,j}
\end{equation}
as $m\rightarrow+\infty$,
where $Z_{n,j}^{\prime}\stackrel{\mathrm{d}}{=}B_{j+\theta/\sigma,
n/\sigma
-j}
Y_{(\theta+n)/\sigma}$ and the random variables $B_{j+\theta/\sigma,
n/\sigma-j}$ and $Y_{(\theta+n)/\sigma}$ are independent. This can be
paralleled with the unconditional limit since it is known that
$K_n/n^\sigma\to Y_{\theta/\sigma}$, almost surely, as $n\to\infty$.
See, for example, \cite{Pit06}, Theorem~3.8.

\begin{remark*} Note that a normalized $\sigma$-stable
process coincides, in distribution, with both a $\operatorname{NGG}(\sigma,0)$ and a
$\operatorname{PD} (\sigma,0)$ process. Hence, it is no surprise that the two limits
\eqref{result} and \eqref{limittwoparameter} are the same, in
distribution, when $\beta=\theta=0$. Another interesting case is
represented by the normalized generalized gamma process with parameter
$(1/2,\beta)$ which yield to the so-called normalized inverse-Gaussian
processes \cite{Lij05}. In particular, for the $\operatorname{NGG}(1/2,\beta)$
process the density $f_{1/2}$ in \eqref{density} is known explicitly
and the previous expression can be simplified to
\begin{eqnarray*}
f_{Z_{n,j}}(z)&=&\frac{\mathrm{e}^{-(\beta/z)^{1/\sigma}}}{\sum
_{l=0}^{n-1}{n-1\choose l}(-1)^{l}\beta^{l/\sigma}\Gamma(j-
{l}/{\sigma};\beta)}
\frac{\Gamma(n)4^{n-1}z^{j/2-1}}{\uppi^{1/2} \Gamma(2n-j) }\\
&&{}\times\sum_{l=0}^{2n-j-1}{2n-j-1\choose l}(-z)^{l/2}\Gamma
\biggl(n-\frac{j-1+l}{2};z\biggr).
\end{eqnarray*}
\end{remark*}

\subsection{Sampling from the limiting random variable}

Since the above described limiting distributions cannot be easily
handled for practical purposes, it is useful to devise a simulation
algorithm. In this respect, one can adapt, similarly to \cite{silvia},
an exact sampling algorithm recently devised by \cite{Dev09} for
random variate generation from polynomially and exponentially tilted
$\sigma$-stable distributions. This will allow to sample the limiting
random variables $Z^{\prime}_{n,j}$ and $Z_{n,j}$ corresponding to the
$\operatorname{PD} (\sigma,\theta)$ and to the $\operatorname{NGG}(\sigma,\beta)$ case, respectively.
Indeed, note that $Z_{n,j}^{\prime}$ is a scale mixture involving a beta
random variable $B_{j+\theta/\sigma, n/\sigma-j}$ and a positive
random variable $Y_{(\theta+n)/\sigma}$. The latter is such that its
transformation $Y_{(\theta+n)/\sigma}^{-1/\sigma}$ admits density
function of the form
%
%
\begin{equation}\label{polynomialtilted}
f_{Y_{(\theta+n)/\sigma}^{-1/\sigma}}(y)=\frac{\Gamma(\theta
+n+1)}{\Gamma({(\theta+n)}/{\sigma}+1)}y^{-\theta-n}f_{\sigma
}(y)\mathbh{1}_{(0,\infty)}(y),
\end{equation}
which is precisely the density function of a polynomially tilted
$\sigma
$-stable distribution. Therefore, random variate generation from
$Z_{n,j}^{\prime}$ can be easily done by independently sampling from a beta
random variable with parameter $(j+\theta/\sigma,n/\sigma-j)$ and from
a random variable with density function \eqref{polynomialtilted} by
means of the algorithm devised in \cite{Dev09}. We refer to \cite
{Fav09} for an alternative sampling algorithm for $Z^{\prime}_{n,j}$ via
augmentation. Similar arguments can be applied
in order to sample from the limit random variable $Z_{n,j}$. Indeed,
observe that $Z_{n,j}$ is characterized by a density function
proportional to
\[
\mathrm{e}^{-(\beta/z)^{1/\sigma}}g_{S_{n,j}}(z)
\]
with $g_{S_{n,j}}$ being the density function of the random variable
$S_{n,j}\stackrel{\mathrm{d}}{=} B_{j,n/\sigma-j}Y_{n/\sigma}$.
Therefore, in order to sample from the distribution of $Z_{n,j}$ one
can apply a simple rejection sampling. In particular, the sampling
scheme would work as follows
\begin{enumerate}[(1)]
\item[(1)] Generate $B \sim B_{j,n/\sigma-j}$.
\item[(2)] Sample $Y \sim Y_{n/\sigma}^{-1/\sigma}$ according to
Devroye's algorithm.\vadjust{\goodbreak}
\item[(3)] Set $S=B Y^{-\sigma}$.
\item[(4)] Sample $U$ from a uniform on the interval $(0,1)$.
\begin{enumerate}[(4.a)]
\item[(4.a)] If $U\le\exp\{-(\beta/S)^{1/\sigma}\}$ set $Z=S$.
\item[(4.b)] If $U>\exp\{-(\beta/S)^{1/\sigma}\}$ restart from (1).
\end{enumerate}
\end{enumerate}

\subsection{Interpretation of asymptotic quantities}
In this final section, we provide a result that gives an interesting
representation of the key random variable $L_{\sigma,n}\stackrel
{\mathrm
{d}}{=}S_{n,j}^{-1/\sigma}$. To this end, we need to provide a
representation for the posterior Laplace transform of the total mass of
the $\sigma$-stable CRM $\tilde\mu_{\sigma,0}$ or, equivalently, of
the unnormalized $\operatorname{NGG}(\sigma,0)$ or $\operatorname{PD} (\sigma, 0)$ processes. Indeed
one has



\begin{prp}\label{poststable}
Let $(X_i)_{i\ge1}$ be a species sampling sequence directed by a
normalized $\sigma$-stable process prior and suppose that the sample
$X_1,\ldots,X_n$ is such that $K_n=j$. Then
%
%
\begin{equation}
\label{eqpostlapl}
\mathbb{E}\bigl[\mathrm{e}^{-\lambda\tilde\mu_{\sigma,0}(\mathbb{X})}
| X_1,\ldots
,X_n\bigr]
=\frac{1}{\Gamma(j)}\int_{\lambda^\sigma}^\infty y^{j-1}\biggl(1-\frac
{\lambda^{1/\sigma}}{y^{1/\sigma}}\biggr)^{n-1} \mathrm{e}^{-y} \,\mathrm{d}y
\end{equation}
for any $\lambda>0$.
\end{prp}

\begin{pf} Set $T_{\sigma,0}\stackrel{\mathrm{d}}{=}\tilde\mu
_{\sigma
,0}(\mathbb{X})$. Since the joint distribution of $(K_n, N_1, \ldots,
N_{K_n})$, also known as exchangeable partition probability function
(see \cite{Pit06}), of a normalized $\sigma$-stable process
coincides with $\mathbb{P}[(K_n, N_1, \ldots, N_{K_n})=(k, n_1,\ldots
,n_k)]=\sigma^{j-1}\Gamma(j) \prod_{i=1}^k(1-\sigma)_{n_i-1}/\Gamma
(n)$, one has
\begin{eqnarray*}
\mathbb{E}[\mathrm{e}^{-\lambda T_{\sigma,0}} | X_1,\ldots,X_n]&=&
\frac{\Gamma(n)}{\sigma^{j-1} \Gamma(j) \prod_{i=1}^k(1-\sigma
)_{n_i-1}}
\frac{1}{\Gamma(n)}\\
&&{}\times\int_0^\infty u^{n-1}\mathrm{e}^{-(\lambda+u)^\sigma} \sigma^j
\prod_{i=1}^k
\frac{\Gamma(n_i-\sigma)}{\Gamma(1-\sigma)} (u+\lambda
)^{-n_i+\sigma
} \,\mathrm{d}u
\end{eqnarray*}
and a simple change of variable $(u+\lambda)^\sigma=y$ yields the
representation in \eqref{eqpostlapl}.
\end{pf}

Proposition~\ref{poststable} allows one to draw an interesting comparison
between unconditional and conditional limits of the number of distinct species.
As we have already highlighted in Section~\ref{section2}, the probability
distribution of the
$\sigma$--diversities for the $\operatorname{NGG}(\sigma,\beta)$ process and the
$\operatorname{PD} (\sigma,\theta)$ process arise
as a power transformation (involving the parameter $\sigma$) of a
suitable tilting of the
probability distribution of $T_{\sigma,0}:=\tilde\mu_{\sigma
,0}(\mathbb{X})$.
We are now in the position to show that
a similar structure carries over when one deals with the conditional
case. Resorting to
the notation set forth in Theorem~\ref{limitresult}, let 
$T_{\sigma,0,K_n}$ to be a random variable whose law coincides with the
probability
distribution of the conditional total mass $T_{\sigma,0}$ of a $\sigma
$-stable process given a sample of size $n$ containing $K_n$ distinct
species. Hence, from
the Laplace transform\vadjust{\goodbreak} \eqref{condition} in the proof of Theorem \ref
{limitresult}
one can easily spot the following identity\vspace*{-2pt}
%
%
\begin{equation}\label{posttotalmass}
L_{\sigma,n}\stackrel{\mathrm{d}}{=} T_{\sigma,0,K_n}.\vadjust{\goodbreak}
\end{equation}
Let now $\mathbb{P}_{\sigma,0}^{(n)}$ and $\mathbb{P}^{(n)}_{\sigma
,\beta}$ be the
conditional probability distributions of, respectively, the
$\sigma$-stable $\tilde{\mu}_{\sigma,0}$ and the generalized gamma
$\tilde{\mu}_{\sigma,\beta}$ processes.
According to Theorem~\ref{limitresult}, the probability distributions
$\mathbb{P}_{\sigma,0}^{(n)}$
and $\mathbb{P}^{(n)}_{\sigma,\beta}$
are mutually absolutely continuous giving rise to the conditional
counterpart of
the identity \eqref{eqexpontilt}, that is,
%
%
\begin{equation}\label{tiltingesponenziale}
\frac{\mathrm{d}\mathbb{P}^{(n)}_{\sigma,\beta} }{\mathrm
{d}\mathbb{P}^{(n)}_{\sigma,0 }}(\mu
)=\frac{\Gamma(j)}{\sum_{l=0}^{n-1}{n-1\choose l}(-1)^{l}\beta
^{l/\sigma
}\Gamma(j-{l}/{\sigma};\beta)}\exp\{-\beta^{1/\sigma}\mu
(\mathbb{X})\}\vspace*{-1pt}
\end{equation}
for any $\sigma\in(0,1)$ and $\beta>0$. In particular, if we denote by
$T_{\sigma,\beta,{K_n}}$ the random
variable whose probability distribution is obtained by exponentially
tilting the
probability distribution of $T_{\sigma,0,K_n}$ as in \eqref
{tiltingesponenziale}, then one
can establish that
%
%
\begin{equation}\label{exptiltngg}
Z_{n,j}\stackrel{\mathrm{d}}{=}(T_{\sigma,\beta,K_n})^{-\sigma}.\vspace*{-1pt}
\end{equation}
In other terms, one can easily verify that the probability distribution
of the limit
random variable $Z_{n,j}$ in \eqref{result} can be also derived by
applying to the
probability distribution of $T_{\sigma,\beta,K_n}$ the same
transformation characterizing the
corresponding unconditional case. In a similar fashion, one can also
derive the
conditional counterpart of the identity \eqref{eqpolytilt} for the
two parameter
Poisson--Dirichlet process. Indeed, according to \cite{Fav09}, Proposition
2,
one can introduce a probability measure $\mathbb{Q}^{(n)}_{\sigma
,\theta
}$ on
$\mathcal{M}_{\mathbb{X}}$ whose Radon--Nikod\'ym derivative with
respect to the dominating
measure $\mathbb{P}_{\sigma,0}^{(n)}$ is given by
%
%
\begin{equation}\label{tiltingpolinomiale}
\frac{\mathrm{d}\mathbb{Q}^{(n)}_{\sigma,\theta} }{\mathrm
{d}\mathbb{P}^{(n)}_{\sigma,0
}}(\mu)=\frac{\Gamma(\theta+n)\Gamma(j)}{\Gamma(n)\Gamma(
{\theta}/{\sigma}+j)}[\mu(\mathbb{X})]^{-\theta}\vspace*{-1pt}
\end{equation}
for any $\sigma\in(0,1)$ and $\theta>-\sigma$ with $\mathbb
{Q}^{(n)}_{\sigma,\theta}$ being the
probability measure of the random measure $\mu^{\ast}_{\sigma,\theta}$
conditional on the sample. In particular, if we denote by
$T^{\prime}_{\sigma
,\theta,K_n}$ the random variable whose probability
distribution is obtained by polynomially tilting the probability
distribution of
$T_{\sigma,0,K_n}$ as in \eqref{tiltingpolinomiale}, then one can
easily verify that
%
%
\begin{equation}\label{exptilttwoparam}
Z^{\prime}_{n,j}\stackrel{\mathrm{d}}{=}(T^{\prime}_{\sigma,\theta
,K_n})^{-\sigma}.\vspace*{-1pt}
\end{equation}
This suggests that the probability distribution of the limiting
random variable $Z_{n,j}^{\prime}$ in \eqref{limittwoparameter} can also be
derived by
applying to the probability distribution of $T'_{\sigma,\theta,K_n}$
the same transformation
characterizing the corresponding unconditional case.\vspace*{-2pt}

\section{Concluding remarks}\vspace*{-2pt}
The identities \eqref{exptiltngg} and \eqref{exptilttwoparam}
represent the
conditional counterparts of the identities \eqref{prima} and~\eqref{seconda},
respectively, given a sample containing $K_n$ distinct species. Hence,
in the same spirit of \cite{Pit03}, Proposition\vadjust{\goodbreak} 13, we have provided
a characterization of the distribution of the limiting random variables
$Z_{n,j}$ and
$Z_{n,j}^{\prime}$ in terms of a power transformation (involving the
parameter~$\sigma$) applied
to a suitable tilting for the conditional distribution of the total
mass of the
$\sigma$-stable process. In particular, the identities \eqref{exptiltngg} and
\eqref{exptilttwoparam} characterize the distribution of the limit random
variables $Z_{n,j}$ and $Z_{n,j}^{\prime}$ via the same transformation characterizing
the unconditional case and applied to an exponential tilting and
polynomial tilting,
respectively, for a scale--mixture distribution involving the beta
distribution and
the $\sigma$-stable distribution. 
To conclude, there is a connection between the prior, and posterior,
total mass of a $\sigma$-stable CRM that we conjecture can be extended
to any
Gibbs-type random probability measure and will be object of future
research.

\section*{Acknowledgements} The authors are grateful to an
Associate Editor and a Referee for valuable remarks and suggestions
that have lead to a substantial improvement in the presentation. This
work is partially supported by MIUR, Grant 2008MK3AFZ, and Regione Piemonte.


%

\printhistory

\end{document}